# Estimation in models driven by fractional Brownian motion

## Corinne Berzin[a] and José R. León[b]

[a] LabSAD – Université Pierre Mendès-France 1251 Avenue centrale, BP 47, 38040 Grenoble cedex 9, France.
E-mail: Corinne.Berzin@upmf-grenoble.fr

[b] Escuela de Matemáticas, Facultad de Ciencias, Universidad Central de Venezuela, Paseo Los Ilustres, Los Chaguaramos, A.P. 47197, Caracas 1041-A, Venezuela. E-mail: jleon@euler.ciens.ucv.ve

**Abstract.** Let $\{b_H(t), t \in \mathbb{R}\}$ be the fractional Brownian motion with parameter $0 < H < 1$. When $1/2 < H$, we consider diffusion equations of the type

$$X(t) = c + \int_0^t \sigma(X(u)) \, db_H(u) + \int_0^t \mu(X(u)) \, du.$$

In different particular models where $\sigma(x) = \sigma$ or $\sigma(x) = \sigma x$ and $\mu(x) = \mu$ or $\mu(x) = \mu x$, we propose a central limit theorem for estimators of $H$ and of $\sigma$ based on regression methods. Then we give tests of the hypothesis on $\sigma$ for these models. We also consider functional estimation on $\sigma(\cdot)$ in the above more general models based in the asymptotic behavior of functionals of the 2nd-order increments of the fBm.

**Résumé.** Soit $\{b_H(t), t \in \mathbb{R}\}$ le mouvement Brownien fractionnaire de paramètre $0 < H < 1$. Lorsque $1/2 < H$, nous considérons des équations de diffusion de la forme

$$X(t) = c + \int_0^t \sigma(X(u)) \, db_H(u) + \int_0^t \mu(X(u)) \, du.$$

Nous proposons dans des modèles particuliers où, $\sigma(x) = \sigma$ ou $\sigma(x) = \sigma x$ et $\mu(x) = \mu$ ou $\mu(x) = \mu x$, un théorème central limite pour des estimateurs de $H$ et de $\sigma$, obtenus par une méthode de régression. Ensuite, pour ces modèles, nous proposons des tests d'hypothèses sur $\sigma$. Enfin, dans les modèles plus généraux ci-dessus nous proposons des estimateurs fonctionnels pour la fonction $\sigma(\cdot)$ dont les propriétés sont obtenues via la convergence de fonctionnelles des accroissements doubles du mBf.

## 1. Introduction

Let $\{b_H(t), t \in \mathbb{R}\}$ be the fractional Brownian motion with Hurst coefficient $H$, $0 < H < 1$. For $H > 1/2$, the integral with respect to fBm can be defined pathwise as the limit of Riemann sums (see [8] and [9]). This allows us to consider, under certain restrictions over $\sigma(\cdot)$ and $\mu(\cdot)$, the "pseudo-diffusion" equations with respect to fBm, that is,

$$X(t) = c + \int_0^t \sigma(X(u))\,\mathrm{d}b_H(u) + \int_0^t \mu(X(u))\,\mathrm{d}u. \tag{1}$$

Our main interest in this work is to provide estimators for the function $\sigma(\cdot)$.

In Section 3.1 we propose simultaneous estimators of $H$ and of $\sigma$ in models such that $\sigma(x) = \sigma$ or $\sigma(x) = \sigma x$ and $\mu(x) = \mu$ or $\mu(x) = \mu x$. Following [4, 7] and [8], we can give an explicit expression for the unique solution to each equation that is a function of the trajectory of $b_H(\cdot)$. Instead of the original process, we use a mollified version and we assume we observe a smoothed by convolution process, defined as $X_\varepsilon(t) = \varphi_\varepsilon * X(t)$. Here $\varepsilon$, which tends to zero, is the smoothing parameter and $\varphi_\varepsilon(\cdot)$ is a convolution kernel such that $\varphi_\varepsilon(\cdot) = \frac{1}{\varepsilon}\varphi(\frac{\cdot}{\varepsilon})$, with $\varphi(\cdot)$ a $C^2$ positive kernel with $L^1$ norm equal to one. Then we observed functionals of the type $\int_0^1 h(X_\varepsilon(u))|\ddot{X}_\varepsilon(u)|^k\,\mathrm{d}u$, with $h(x) = 1/|x|^k$ in the case of linear $\sigma(\cdot)$ and $h(x) = 1$ in the case of constant $\sigma(\cdot)$. Such an observation could seem unusual but note that in case where $\varphi(\cdot) = \mathbf{1}_{[-1,0]} * \mathbf{1}_{[-1,0]}(\cdot)$, then $\varepsilon^2 \ddot{X}_\varepsilon(u) = X(u+2\varepsilon) - 2X(u+\varepsilon) + X(u)$. Although $\varphi(\cdot)$ is not a $C^2$ function, all results obtained in this paper can be rewritten with this particular approximation. So our method can be related with variation methods and consists in obtaining some least squares estimators for $H$ and $\sigma$ in certain regression models. This method can be compared to that of [3], where the model was the fBm i.e. $\sigma(\cdot) \equiv 1$ and $\mu(\cdot) \equiv 0$ and the purpose was to estimate $H$. Indeed, we prove that the asymptotic behavior of such estimators, that is, $(\widehat{H}_k - H)/\sqrt{\varepsilon}$ and $(\widehat{\sigma}_k - \sigma)/(\sqrt{\varepsilon}\log(\varepsilon))$ are both equivalent to those of certain non-linear functionals of the Gaussian process $b_H^\varepsilon(\cdot) = \varphi_\varepsilon * b_H(\cdot)$. As in [3], we show that they satisfy a central limit theorem using the method of moments, via the diagram formula. It is interesting to note that the rates of convergence of such estimators are not the same. Furthermore the asymptotic law of $(\widehat{H}_k, \widehat{\sigma}_k)$ is a degenerated Gaussian. Hence, we could not provide simultaneous confidence intervals for $H$ and $\sigma$. Finally, we proved that the best estimators for $H$ and $\sigma$ in the sense of minimal variance are obtained for $k = 2$.

In Section 3.2, we get back to more general models of the form (1) and our goal is to provide functional estimation of $\sigma(\cdot)$ as in [2]. Indeed, in [2] we considered the case where $\mu(\cdot) \equiv 0$ and we proved that, if $N_\varepsilon^X(x)$ denotes the number of times the regularized process $X_\varepsilon(\cdot)$ crosses level $x$, before time 1, then for $1/2 < H < 1$ and any continuous function $h$,

$$\sqrt{\frac{\pi}{2}}\frac{\varepsilon^{(1-H)}}{\tilde{\sigma}_{2H}}\int_{-\infty}^\infty h(x) N_\varepsilon^X(x)\,\mathrm{d}x \xrightarrow{\text{a.s.}} \int_0^1 h(X(u))\sigma(X(u))\,\mathrm{d}u. \tag{2}$$

Furthermore we got the following result about the rates of convergence proving that there exists a Brownian motion $\widehat{W}(\cdot)$ independent of $b_H(\cdot)$ and a constant $\sigma_{g_1}$ such that for $1/2 < H < 3/4$,

$$\frac{1}{\sqrt{\varepsilon}}\left[\sqrt{\frac{\pi}{2}}\frac{\varepsilon^{(1-H)}}{\tilde{\sigma}_{2H}}\int_{-\infty}^\infty h(x) N_\varepsilon^X(x)\,\mathrm{d}x - \int_0^1 h(X(u))\sigma(X(u))\,\mathrm{d}u\right]$$
$$\underset{\varepsilon \to 0}{\Longrightarrow} \sigma_{g_1}\int_0^1 h(X(u))\sigma(X(u))\,\mathrm{d}\widehat{W}(u), \tag{3}$$

under some assumptions on the function $h$.

The proofs of these two last convergences are based on the fact that, on the one hand, when $\mu(\cdot) \equiv 0$ and because fBm has quadratic variation when $H > \frac{1}{2}$ and $\sigma(\cdot) \in C^1(\mathbb{R})$, the solution to Eq. (1) is given by $X(t) = K(b_H(t))$ where $K$ is the solution to the ordinary differential equation $\dot{K}(t) = \sigma(K(t))$ with $K(0) = c$ (see [8]). On the other hand, by the Banach formula, we have $\int_{-\infty}^\infty h(x) N_\varepsilon^X(x)\,\mathrm{d}x = \int_0^1 h(X_\varepsilon(u))|\dot{X}_\varepsilon(u)|\,\mathrm{d}u$ and since $\dot{X}_\varepsilon(u) \simeq \dot{K}(b_H^\varepsilon(u))\dot{b}_H^\varepsilon(u) = \sigma(K(b_H^\varepsilon(u)))\dot{b}_H^\varepsilon(u)$, we needed to look for the asymptotic behavior of a particular non-linear functional of the regularized fBm $b_H^\varepsilon(\cdot)$ and Theorems 1.1 and 3.4 of [2] gave the result.

We need to recall these two convergence results because they can serve as a motivation to the statement of Theorems 3.8 and 3.9. Indeed in these two last theorems we will use the same type of approximations and convergence. Also convergences in (2) and (3) can serve as a motivation to the main interest in the present paper, that is to work with the second derivative of the smoothed process instead of the first one. Indeed in the second convergence result (3), we obtained the restrictive condition $H < 3/4$ due to the fact that we used the first derivative of $X_\varepsilon(\cdot)$.

Thus in this work, having in mind to reach all the range of $H < 1$, we considered the case where $\mu(\cdot)$ is not necessarily null and instead of considering functionals of $|\dot{X}_\varepsilon(\cdot)|$, we worked with functionals of $|\ddot{X}_\varepsilon(\cdot)|^k$. This approach allowed us to provide functional estimation of $\sigma(\cdot)$ as in (2) and to exhibit the rate of convergence as in (3) for any value of $H$ in $]1/2, 1[$, using a generalization of the two last convergence results when $\mu(\cdot) \equiv 0$ and then applying the Girsanov theorem (see [5] and [6]).

We observe that the limit convergence in (2) is $\int_0^1 h(X(u))\sigma(X(u))\,\mathrm{d}u$ and will become $\int_0^1 \sigma^k(X(u))\,\mathrm{d}u$ in that work in cases where $h(\cdot) \equiv 1$. Thus if we get back to the last four models of Section 3.1 and if we take into account the form of $\sigma(\cdot)$, that is $\sigma(x) = \sigma$ or $\sigma(x) = \sigma x$, the limit integral is now a function of $\sigma$. Thus in the second part of Section 3.1 we propose estimators for $\sigma^k$ when $H$ is known. This supplementary information about $H$ leads us to estimators of $\sigma$ performing more than those of the first part of Section 3.1 because the rate of convergence will be $1/\sqrt{\varepsilon}$ instead of $1/(\sqrt{\varepsilon}\log(\varepsilon))$ as before.

Finally in Section 3.3, a result similar to the one of part two in Section 3.1 can be obtained under contiguous alternatives for $\sigma$ and provides a test of the hypothesis for such a coefficient.

We work with the techniques described in the last sections but without using the Girsanov theorem that is more tricky to use since under the alternative hypothesis $\sigma$ depends on $\varepsilon$.

## 2. Hypothesis and notation

Let $\{b_H(t), t \in \mathbb{R}\}$ be a fractional Brownian motion (fBm) with parameter $0 < H < 1$ (see for instance [10]), i.e. $b_H(\cdot)$ is a centered Gaussian process with the covariance function:

$$\mathbb{E}[b_H(t)b_H(s)] = \frac{1}{2}v_{2H}^2[|t|^{2H} + |s|^{2H} - |t-s|^{2H}],$$

with $v_{2H}^2 := [\Gamma(2H+1)\sin(\pi H)]^{-1}$. We define, for a $C^2$ density $\varphi$ with compact support included in $[-1, 1]$, for each $t \geq 0$ and $\varepsilon > 0$ the regularized processes:

$$b_H^\varepsilon(t) := \frac{1}{\varepsilon}\int_{-\infty}^{\infty} \varphi\left(\frac{t-x}{\varepsilon}\right) b_H(x)\,\mathrm{d}x \quad \text{and} \quad Z_\varepsilon(t) := \frac{\varepsilon^{(2-H)}\ddot{b}_H^\varepsilon(t)}{\sigma_{2H}},$$

with

$$\sigma_{2H}^2 := \mathbb{V}[\varepsilon^{(2-H)}\ddot{b}_H^\varepsilon(t)] = \frac{1}{2\pi}\int_{-\infty}^{+\infty} |x|^{(3-2H)}|\hat{\varphi}(-x)|^2\,\mathrm{d}x.$$

We shall use Hermite polynomials, defined by

$$\mathrm{e}^{(tx-t^2/2)} = \sum_{n=0}^{\infty} \frac{H_n(x)t^n}{n!}.$$

They form an orthogonal system for the standard Gaussian measure $\phi(x)\,\mathrm{d}x$ and, if $h \in L^2(\phi(x)\,\mathrm{d}x)$, then $h(x) = \sum_{n=0}^{\infty} \hat{h}_n H_n(x)$ and $\|h\|_{2,\phi}^2 := \sum_{n=0}^{\infty} \hat{h}_n^2 n!$.

Let $g$ be a function in $L^2(\phi(x)\,\mathrm{d}x)$ such that $g(x) = \sum_{n=1}^{\infty} \hat{g}_n H_n(x)$, with $\|g\|_{2,\phi}^2 = \sum_{n=1}^{\infty} \hat{g}_n^2 n! < +\infty$.

The symbol "$\Rightarrow$" will mean weak convergence of measures.

At this step of the paper it will be helpful to state several theorems obtained in [3] in the aim to enlighten the notations and to make this paper more independent of it.

We proved the two following theorems.

**Theorem 2.1.** *For all $0 < H < 1$ and $k \in \mathbb{N}^*$,*

$$\int_0^1 [Z_\varepsilon(u)]^k \, du \xrightarrow[\varepsilon \to 0]{a.s.} \mathbb{E}[N]^k,$$

*where $N$ will denote a standard Gaussian random variable.*

**Remark.** *This theorem implies that $Z_\varepsilon(\cdot) \Rightarrow N$ when $\varepsilon$ goes to zero, the random variable $Z_\varepsilon(\cdot)$ is considered as a variable on $([0,1], \lambda)$ where $\lambda$ is the Lebesgue measure. This last convergence implies that for all $0 < H < 1$ and real $k > 0$, almost surely, $\int_0^1 |Z_\varepsilon(u)|^k \, du \to \mathbb{E}[|N|^k]$.*

For $\varepsilon > 0$, define

$$S_g(\varepsilon) := \varepsilon^{-1/2} \int_0^1 g(Z_\varepsilon(u)) \, du.$$

**Theorem 2.2.** *For all $0 < H < 1$,*

$$S_g(\varepsilon \cdot) \underset{\varepsilon \to 0}{\Longrightarrow} X(\cdot),$$

*where the above convergence is in the sense of finite-dimensional distributions and $X(\cdot)$ is a cylindrical centered Gaussian process with covariance $\rho_g(b,c) := \mathbb{E}[X(b)X(c)]$, where for $b, c > 0$,*

$$\rho_g(b,c) := \frac{1}{\sqrt{bc}} \sum_{n=1}^{\infty} \hat{g}_n^2 n! \int_{-\infty}^{+\infty} \rho_H^n(x,b,c) \, dx,$$

*and for $x \in \mathbb{R}$,*

$$\rho_H(x,b,c) := \mathbb{E}[Z_{\varepsilon b}(\varepsilon x + u) Z_{\varepsilon c}(u)] = \frac{(bc)^{(2-H)}}{2\pi \sigma_{2H}^2} \int_{-\infty}^{\infty} |y|^{(3-2H)} e^{ixy} \widehat{\varphi}(-by) \widehat{\varphi}(cy) \, dy.$$

Note that for fixed $c > 0$, $\rho_H(x,c,c) = \rho_H(x/c)$, where for $x \in \mathbb{R}$,

$$\rho_H(x) := \mathbb{E}[Z_\varepsilon(\varepsilon x + u) Z_\varepsilon(u)] = \frac{1}{2\pi \sigma_{2H}^2} \int_{-\infty}^{\infty} |y|^{(3-2H)} e^{ixy} |\widehat{\varphi}(-y)|^2 \, dy,$$

so that, $\rho_g(c,c) = \sigma_g^2 := \sum_{n=1}^{\infty} \hat{g}_n^2 n! \int_{-\infty}^{+\infty} \rho_H^n(x) \, dx$. Therefore we furthermore got the following remark that was also shown in Corollary 3.2(i) in [2].

**Remark.** *If $c > 0$ is fixed, $S_g(\varepsilon c) \Rightarrow \sigma_g N$.*

For all $m \in \mathbb{N}^*$, for all $c_1 > 0$, $c_2 > 0$, ..., $c_m > 0$ and for all $d_1, d_2, \ldots, d_m \in \mathbb{R}$, we will denote

$$\sigma_{g,m}^2(\mathbf{c}, \mathbf{d}) := \sum_{i=1}^{m} \sum_{j=1}^{m} d_i d_j \rho_g(c_i, c_j) = \mathbb{E}\left[\sum_{i=1}^{m} d_i X(c_i)\right]^2.$$

We also proved in [3] the following lemma.

**Lemma 2.1.** $\lim_{\varepsilon \to 0} \mathbb{E}[\sum_{i=1}^{m} d_i S_g(\varepsilon c_i)]^2 = \sigma_{g,m}^2(\mathbf{c}, \mathbf{d})$.

Throughout the paper, $\mathbf{C}$ (resp. $\mathbf{C}(\omega)$) shall stand for a generic constant (resp. for a generic constant that depends on $\omega$ living in the space of the trajectories), whose value may change during a proof, and $\log(\cdot)$ for the Naperian logarithm.

For $k \geq 1$, we shall note $\|N\|_k^k := \mathbb{E}[|N|^k]$, and if $C(\cdot)$ is a measurable function we shall note $\|C(\cdot)\|_k^k := \int_0^1 |C(u)|^k \,\mathrm{d}u$, by $\|C(\cdot)\|_{k,\varepsilon}^k := \int_0^\varepsilon |C(u)|^k \,\mathrm{d}u$ and by $\|C(\cdot)\|_{k,\varepsilon^c}^k := \int_\varepsilon^1 |C(u)|^k \,\mathrm{d}u$, for $\varepsilon > 0$.

## 3. Results

### 3.1. Simultaneous estimators of $H$ and of $\sigma$

As mentioned in the Introduction, we are interested in providing simultaneous estimators of $H$ and $\sigma$ in the four following models. For $H > 1/2$ and $t \geq 0$

$$\mathrm{d}X(t) = \sigma \,\mathrm{d}b_H(t) + \mu \,\mathrm{d}t, \tag{4}$$

$$\mathrm{d}X(t) = \sigma \,\mathrm{d}b_H(t) + \mu X(t) \,\mathrm{d}t, \tag{5}$$

$$\mathrm{d}X(t) = \sigma X(t) \,\mathrm{d}b_H(t) + \mu X(t) \,\mathrm{d}t, \tag{6}$$

$$\mathrm{d}X(t) = \sigma X(t) \,\mathrm{d}b_H(t) + \mu \,\mathrm{d}t, \tag{7}$$

with $X(0) = c$.

The solution to these equations are, respectively:

(4) $X(t) = \sigma b_H(t) + \mu t + c$ (see [8]),
(5) $X(t) = \sigma b_H(t) + \exp(\mu t)[\sigma \mu(\int_0^t b_H(s) \exp(-\mu s) \,\mathrm{d}s) + c]$,
(6) $X(t) = c \exp(\mu t + \sigma b_H(t))$ (see [4] and [7]),
(7) $X(t) = \exp(\sigma b_H(t))(c + \mu \int_0^t \exp(-\sigma b_H(s)) \,\mathrm{d}s)$.

We consider the problem of estimating simultaneously $H$ and $\sigma > 0$. Suppose we observe instead of $X(t)$ a smoothed by convolution process $X_\varepsilon(t) = \varphi_\varepsilon * X(t)$, where $\varphi_\varepsilon(\cdot) = 1/\varepsilon \varphi(\cdot/\varepsilon)$, with $\varphi(\cdot)$ as before and where we have extended $X(\cdot)$ by means of $X(t) = c$, if $t < 0$.

For model (7) we will make the additional hypothesis that $\mu$ and $c$ have the same sign with $\mu$ possibly null.

From now on, we shall note for each $t \geq 0$ and $\varepsilon > 0$,

$$Z_\varepsilon^X(t) := \begin{cases} \dfrac{\varepsilon^{(2-H)} \ddot{X}_\varepsilon(t)}{\sigma 2H} & \text{for the first two models,} \\ \dfrac{\varepsilon^{(2-H)} \ddot{X}_\varepsilon(t)}{\sigma 2H X_\varepsilon(t)} & \text{for the other two.} \end{cases} \tag{8}$$

For $k \geq 1$, let us denote

$$A_k(\varepsilon) := \frac{1}{\sigma^k \|N\|_k^k} \left( \int_0^1 |Z_\varepsilon^X(u)|^k \,\mathrm{d}u \right) - 1. \tag{9}$$

The remark following Theorem 2.1 allows us to state the following theorem.

**Theorem 3.1.** *For $k \geq 1$,*

(1) $A_k(\varepsilon) \xrightarrow[\varepsilon \to 0]{a.s.} 0$.
(2) *Furthermore*

$$\frac{1}{\sqrt{\varepsilon}} A_k(\varepsilon) = S_{g_k}(\varepsilon) + \mathrm{o}_{a.s.}(1), \quad \text{where}$$

$$g_k(x) := \frac{|x|^k}{\|N\|_k^k} - 1. \tag{10}$$

At this step, we can propose estimators of $H$ and $\sigma$, by observing $X_\varepsilon(u)$ at several scales of the parameter $\varepsilon$, i.e. $h_i = \varepsilon c_i$, $c_i > 0$, $i = 1, \ldots, \ell$. In this aim, let us define

$$M_k(\varepsilon) := \begin{cases} \int_0^1 |\ddot{X}_\varepsilon(u)|^k \, du & \text{for the first two models,} \\ \int_0^1 |\frac{\ddot{X}_\varepsilon(u)}{X_\varepsilon(u)}|^k \, du & \text{for the other two.} \end{cases} \tag{11}$$

Using assertion (1) of Theorem 3.1, we get

$$\frac{\varepsilon^{k(2-H)} M_k(\varepsilon)}{\sigma_{2H}^k \sigma^k \|N\|_k^k} \xrightarrow[\varepsilon \to 0]{\text{a.s.}} 1,$$

from which we obtain

$$\log(M_k(\varepsilon)) = k(H-2)\log(\varepsilon) + \log(\sigma_{2H}^k \sigma^k \|N\|_k^k) + o_{\text{a.s.}}(1). \tag{12}$$

The following regression model can be written, for each scale $h_i$:

$$Y_i = a_k X_i + b_k + \xi_i, \quad i = 1, \ldots, \ell,$$

where $a_k := k(H-2)$, $b_k := \log(\sigma_{2H}^k \sigma^k \|N\|_k^k)$ and for $i = 1, \ldots, \ell$, $Y_i := \log(M_k(h_i))$, $X_i := \log(h_i)$. Hence, the least squares estimators $\widehat{H}_k$ of $H$ and $\widehat{B}_k$ of $b_k$ are defined as

$$k(\widehat{H}_k - 2) := \sum_{i=1}^\ell z_i \log(M_k(h_i)) \tag{13}$$

and

$$\widehat{B}_k := \frac{1}{\ell} \sum_{i=1}^\ell \log(M_k(h_i)) - k(\widehat{H}_k - 2)\frac{1}{\ell}\sum_{i=1}^\ell \log(h_i), \tag{14}$$

where

$$z_i := \frac{y_i}{\sum_{i=1}^\ell y_i^2} \quad \text{and} \quad y_i := \log(c_i) - \frac{1}{\ell}\sum_{i=1}^\ell \log(c_i).$$

Note the following property

$$\sum_{i=1}^\ell z_i = 0 \quad \text{and} \quad \sum_{i=1}^\ell z_i y_i = 1. \tag{15}$$

Then we propose

$$\widehat{\sigma_{2H}^k} := \sigma_{2\widehat{H}_k}^k,$$

as estimator of $\sigma_{2H}^k$ and

$$\widehat{\sigma^k} := \frac{\exp(\widehat{B}_k)}{\widehat{\sigma_{2H}^k} \|N\|_k^k}, \tag{16}$$

as estimator of $\sigma^k$. Finally, we propose $\widehat{\sigma}_k$ as estimator of $\sigma$ defined by

$$\widehat{\sigma}_k := (\widehat{\sigma^k})^{1/k}. \tag{17}$$

Theorem 3.1 and Theorem 2.2 imply the following theorem for all the range of $H$ belonging to $]1/2, 1[$.

**Theorem 3.2.** *For $k \geq 1$,*

(1) *$\widehat{H}_k$ is a strongly consistent estimator of $H$ and*

$$\frac{1}{\sqrt{\varepsilon}}(\widehat{H}_k - H) \underset{\varepsilon \to 0}{\Longrightarrow} \mathcal{N}\left(0, \sigma_{g_k,\ell}^2\left(\mathbf{c}, \sqrt{\mathbf{c}}\left(\frac{\mathbf{z}}{k}\right)\right)\right),$$

*where $g_k(\cdot)$ is defined by (10) and*

$$g_k(x) = \sum_{n=1}^{\infty} \hat{g}_{2n,k} H_{2n}(x), \quad \text{with } \hat{g}_{2n,k} = \frac{1}{(2n)!} \prod_{i=0}^{n-1}(k-2i).$$

(2) *$\widehat{\sigma}_k$ is a weakly consistent estimator of $\sigma$ and*

$$\frac{1}{\sqrt{\varepsilon}\log(\varepsilon)}(\widehat{\sigma}_k - \sigma) \underset{\varepsilon \to 0}{\Longrightarrow} \mathcal{N}\left(0, \sigma^2 \sigma_{g_k,\ell}^2\left(\mathbf{c}, \sqrt{\mathbf{c}}\left(\frac{\mathbf{z}}{k}\right)\right)\right).$$

**Remark.** *As in [3], the variance $\sigma_{g_k,\ell}^2(\mathbf{c}, \sqrt{\mathbf{c}}(\mathbf{z}/k))$ is minimal for $k=2$ and then the best estimators for $H$ and $\sigma$ in the sense of minimal variance are obtained for $k=2$.*

Now let us suppose that $H$, $\frac{1}{2} < H < 1$, is known. Theorem 3.1 also provides estimators for $\sigma$. Indeed, if for $k \geq 1$ we set,

$$\tilde{\sigma}_k := \frac{\|Z_\varepsilon^X(\cdot)\|_k}{\|N\|_k} \quad \text{(see (8) for definition of } Z_\varepsilon^X(\cdot)\text{),}$$

then Theorem 3.1 and the remark following Theorem 2.2 imply the following theorem.

**Theorem 3.3.** *For $k \geq 1$ and if $H$ is known with $\frac{1}{2} < H < 1$ then*

(1) *$\tilde{\sigma}_k$ is a strongly consistent estimator of $\sigma$ and*
(2)

$$\frac{1}{\sqrt{\varepsilon}}(\tilde{\sigma}_k - \sigma) \underset{\varepsilon \to 0}{\Longrightarrow} \mathcal{N}\left(0, \frac{\sigma^2}{k^2}\sigma_{g_k}^2\right),$$

*where $g_k(\cdot)$ is defined by (10).*

**Remark 1.** *Note that the rate of convergence in assertion (2) is $1/\sqrt{\varepsilon}$ instead of $1/(\sqrt{\varepsilon}\log(\varepsilon))$ as in assertion (2) in Theorem 3.2. This is due to the fact that here $H$ is known.*

**Remark 2.** *The variance $\sigma_{g_k}^2/k^2$ is minimal for $k=2$ and then the best estimator for $\sigma$ in the sense of minimal variance is obtained for $k=2$.*

**Proof of Theorem 3.1.**
(1) We need the following lemma for which proof is given in the Appendix.

**Lemma 3.1.** *For $0 \leq t \leq 1$,*

$$\ddot{X}_\varepsilon(t) = \begin{cases} \sigma \ddot{b}_H^\varepsilon(t) + a_\varepsilon(t) & \text{for the first two models,} \\ \sigma X_\varepsilon(t)\ddot{b}_H^\varepsilon(t) + X_\varepsilon(t)a_\varepsilon(t) & \text{for the other two,} \end{cases}$$

*with*

$$|a_\varepsilon(t)| \leq \mathbf{C}(\omega)(\varepsilon^{(H-2-\delta)}\mathbf{1}_{\{0 \leq t \leq \varepsilon\}} + \varepsilon^{(2H-2-\delta)}\mathbf{1}_{\{\varepsilon \leq t \leq 1\}}), \quad \text{for any } \delta > 0.$$

**Remark.** *Indeed, for the first model one has,*

$$|a_\varepsilon(t)| \leq \mathbf{C}(\omega)\varepsilon^{(H-2-\delta)}\mathbf{1}_{\{0 \leq t \leq \varepsilon\}}$$

*and for the second one,*

$$|a_\varepsilon(t)| \leq \mathbf{C}(\omega)(\varepsilon^{(H-2-\delta)}\mathbf{1}_{\{0 \leq t \leq \varepsilon\}} + \varepsilon^{(H-1-\delta)}\mathbf{1}_{\{\varepsilon \leq t \leq 1\}}).$$

We have to prove that almost surely, $\|Z_\varepsilon^X(\cdot)\|_k$ converges to $\sigma\|N\|_k$ when $\varepsilon$ goes to zero. For this, we write $\|Z_\varepsilon^X(\cdot)\|_k$ as

$$\|Z_\varepsilon^X(\cdot)\|_k = \|\sigma Z_\varepsilon(\cdot)\|_k + \|Z_\varepsilon^X(\cdot)\|_k - \|\sigma Z_\varepsilon(\cdot)\|_k.$$

By the remark following Theorem 2.1, we know that $\|Z_\varepsilon(\cdot)\|_k$ converges almost surely to $\|N\|_k$ when $\varepsilon$ goes to zero. Thus, we have just to prove that $|\|Z_\varepsilon^X(\cdot)\|_k - \|\sigma Z_\varepsilon(\cdot)\|_k|$ converges almost surely to zero with $\varepsilon$. By Lemma 3.1 and using Minkowski's inequality, one has

$$|\|Z_\varepsilon^X(\cdot)\|_k - \|\sigma Z_\varepsilon(\cdot)\|_k| \leq \|Z_\varepsilon^X(\cdot) - \sigma Z_\varepsilon(\cdot)\|_k$$
$$= \frac{\varepsilon^{(2-H)}}{\sigma_{2H}}\|a_\varepsilon(\cdot)\|_k \leq \mathbf{C}(\omega)(\varepsilon^{(1/k-\delta)} + \varepsilon^{(H-\delta)}).$$

Choosing $\delta$ small enough, i.e. $0 < \delta < \inf(\frac{1}{k}, H)$, we proved that the last term in the above inequality tends almost surely to zero with $\varepsilon$ and assertion (1) follows.

(2) We write

$$\frac{1}{\sqrt{\varepsilon}}A_k(\varepsilon) = S_{g_k}(\varepsilon) + \frac{1}{\sqrt{\varepsilon}\|N\|_k^k\sigma^k}(\|Z_\varepsilon^X(\cdot)\|_k^k - \|\sigma Z_\varepsilon(\cdot)\|_k^k).$$

Let us prove now that $(\|Z_\varepsilon^X(\cdot)\|_k^k - \|\sigma Z_\varepsilon(\cdot)\|_k^k) = o_{\mathrm{a.s.}}(\sqrt{\varepsilon})$. Using the bound

$$||x+y|^k - |x|^k| \leq 2^{(k-1)}k|y|(|x|^{(k-1)} + |y|^{(k-1)}), \quad \text{for } k \geq 1,$$

and Hölder's inequality, we obtain

$$|\|f\|_k^k - \|g\|_k^k| \leq \||f(\cdot)|^k - |g(\cdot)|^k\|_1$$
$$\leq 2^{(k-1)}k\|f - g\|_k[\|g\|_k^{(k-1)} + \|f - g\|_k^{(k-1)}]. \tag{18}$$

Let us apply this inequality to $f(\cdot) := Z_\varepsilon^X(\cdot)$ and to $g(\cdot) := \sigma Z_\varepsilon(\cdot)$, successively with the norm $\|\cdot\|_{k,\varepsilon}$ and the norm $\|\cdot\|_{k,\varepsilon^c}$.

On the one hand, applying Lemma 3.1, one obtains that $\|Z_\varepsilon^X(\cdot) - \sigma Z_\varepsilon(\cdot)\|_{k,\varepsilon} \leq \mathbf{C}(\omega)\varepsilon^{(1/k-\delta)}$ and that $\|Z_\varepsilon^X(\cdot) - \sigma Z_\varepsilon(\cdot)\|_{k,\varepsilon^c} \leq \mathbf{C}(\omega)\varepsilon^{(H-\delta)}$.

On the other hand, the trajectories of $b_H(\cdot)$ are $(H-\delta)$-Hölder continuous, in other words for any $\delta > 0$

$$|b_H(u+\varepsilon) - b_H(u)| \leq \mathbf{C}(\omega)\varepsilon^{(H-\delta)}. \tag{19}$$

Using this fact, we get

$$|\ddot{b}_H^\varepsilon(u)| = \left|\frac{1}{\varepsilon^2}\int_{-\infty}^\infty \ddot{\varphi}(v)(b_H(u-\varepsilon v) - b_H(u))\,\mathrm{d}v\right| \leq \mathbf{C}(\omega)\varepsilon^{(H-2-\delta)}. \tag{20}$$

We deduce $\|Z_\varepsilon(\cdot)\|_{k,\varepsilon} \leq \mathbf{C}(\omega)\varepsilon^{(1/k-\delta)}$ and $\|Z_\varepsilon(\cdot)\|_{k,\varepsilon^c} \leq \mathbf{C}(\omega)\varepsilon^{-\delta}$.

Finally, taking $\delta$ small enough, i.e. $0 < \delta < \frac{1}{k}(H - \frac{1}{2})$, we proved that

$$|\|Z_\varepsilon^X(\cdot)\|_k^k - \|\sigma Z_\varepsilon(\cdot)\|_k^k| \leq \mathbf{C}(\omega)(\varepsilon^{(1-\delta k)} + \varepsilon^{(H-\delta k)} + \varepsilon^{k(H-\delta)})$$
$$\leq \mathbf{C}(\omega)\varepsilon^{(H-\delta k)} = o_{\text{a.s.}}(\sqrt{\varepsilon}),$$

and assertion (2) follows. $\square$

**Proof of Theorem 3.2.**

(1) By using (13) and (12) we obtain

$$k(\widehat{H}_k - 2) = \sum_{i=1}^{\ell} z_i[k(H-2)\log(\varepsilon c_i) + b_k] + o_{\text{a.s.}}(1),$$

and property (15) gives

$$k(\widehat{H}_k - 2) = k(H - 2) + o_{\text{a.s.}}(1).$$

We proved that $\widehat{H}_k$ is a strongly consistent estimator of $H$.

Now by using definitions of $A_k(\varepsilon)$ and of $M_k(\varepsilon)$ (see (9) and (11)), one obtains

$$A_k(\varepsilon) = \varepsilon^{k(2-H)} \exp(-b_k) M_k(\varepsilon) - 1.$$

With this definition and using a Taylor expansion for the logarithm function one has

$$\log(M_k(\varepsilon)) = \log(\varepsilon^{k(H-2)} \exp(b_k)) + \log(1 + A_k(\varepsilon))$$
$$= k(H-2)\log(\varepsilon) + b_k + A_k(\varepsilon) + A_k^2(\varepsilon)\left[-\frac{1}{2} + \varepsilon_1(A_k(\varepsilon))\right]. \quad (21)$$

Let us see that

$$A_k^2(\varepsilon)\left[-\frac{1}{2} + \varepsilon_1(A_k(\varepsilon))\right] = o_P(\sqrt{\varepsilon}). \quad (22)$$

By assertion (2) of Theorem 3.1 we know that

$$\frac{1}{\sqrt{\varepsilon}} A_k(\varepsilon) = S_{g_k}(\varepsilon) + o_{\text{a.s.}}(1), \quad (23)$$

where $g_k(\cdot)$ is defined by (10), and by Lemma 2.1,

$$\mathbb{E}[S_{g_k}^2(\varepsilon)] = O(1), \quad (24)$$

so $\frac{1}{\sqrt{\varepsilon}} A_k^2(\varepsilon) = o_P(1)$ and then (22) is proved.

By using (21)–(23) we obtain

$$\log(M_k(\varepsilon)) = k(H-2)\log(\varepsilon) + b_k + \sqrt{\varepsilon} S_{g_k}(\varepsilon) + o_P(\sqrt{\varepsilon}). \quad (25)$$

Thus (13), (25) and property (15) entail that

$$k(\widehat{H}_k - 2) = k(H - 2) + \sum_{i=1}^{\ell} z_i \sqrt{\varepsilon c_i} S_{g_k}(\varepsilon c_i) + o_P(\sqrt{\varepsilon}).$$

Then

$$\frac{(\widehat{H}_k - H)}{\sqrt{\varepsilon}} = \frac{1}{k} \sum_{i=1}^{\ell} z_i \sqrt{c_i} S_{g_k}(\varepsilon c_i) + o_P(1).$$

Theorem 2.2 gives the required result (the computation of the coefficients in the Hermite expansion of function $g_k(\cdot)$ is explicitly made in the proof of Corollary 3.2 of [3]).

(2) Let us see that $\widehat{B}_k$ is a weakly consistent estimator of $b_k$.

By using (14) and (25), one has

$$\widehat{B}_k - b_k = \frac{k}{\ell}(H - \widehat{H}_k)\sum_{i=1}^{\ell}\log(h_i) + \sqrt{\varepsilon}\frac{1}{\ell}\sum_{i=1}^{\ell}\sqrt{c_i}S_{g_k}(\varepsilon c_i) + \sqrt{\varepsilon}o_P(1).$$

Thus

$$\widehat{B}_k - b_k = k\log(\varepsilon)(H - \widehat{H}_k) + \frac{k}{\ell}(H - \widehat{H}_k)\sum_{i=1}^{\ell}\log(c_i) + \sqrt{\varepsilon}\frac{1}{\ell}\sum_{i=1}^{\ell}\sqrt{c_i}S_{g_k}(\varepsilon c_i) + \sqrt{\varepsilon}o_P(1).$$

Using assertion (1) of Theorem 3.2 and (24), we obtain

$$\frac{[\widehat{B}_k - b_k]}{\sqrt{\varepsilon}\log(\varepsilon)} = k\left(\frac{H - \widehat{H}_k}{\sqrt{\varepsilon}}\right) + o_P(1), \tag{26}$$

and then using again assertion (1) of Theorem 3.2 we proved that $\widehat{B}_k$ is a weakly consistent estimator of $b_k$.

Now using a Taylor expansion of order one for the exponential function, equality (26), the fact that $\widehat{B}_k$ is a weakly consistent estimator of $b_k$ and assertion (1) of Theorem 3.2, we finally get

$$\frac{[\exp(\widehat{B}_k) - \exp(b_k)]}{\sqrt{\varepsilon}\log(\varepsilon)} = k\exp(b_k)\left(\frac{H - \widehat{H}_k}{\sqrt{\varepsilon}}\right) + o_P(1). \tag{27}$$

Thus if we get back to the definition of $\widehat{\sigma^k}$ (see (16)) and if we use the last equality (27) we get

$$\frac{(\widehat{\sigma^k} - \sigma^k)}{\sqrt{\varepsilon}\log(\varepsilon)} = \frac{\sigma^k \exp(-b_k)}{\sqrt{\varepsilon}\log(\varepsilon)}\left\{[\exp(\widehat{B}_k) - \exp(b_k)] + \exp(\widehat{B}_k)\left[\frac{1}{\sigma^k_{2\widehat{H}_k}} - \frac{1}{\sigma^k_{2H}}\right]\sigma^k_{2H}\right\}$$

$$= k\sigma^k\left(\frac{H - \widehat{H}_k}{\sqrt{\varepsilon}}\right) + \left(\frac{\sigma}{\sigma_{2\widehat{H}_k}}\right)^k \exp(-b_k)\exp(\widehat{B}_k)\frac{[\sigma^k_{2H} - \sigma^k_{2\widehat{H}_k}]}{\sqrt{\varepsilon}\log(\varepsilon)} + o_P(1).$$

At this step of the proof we are going to show that

$$\frac{(\widehat{\sigma^k} - \sigma^k)}{\sqrt{\varepsilon}\log(\varepsilon)} = k\sigma^k\left(\frac{H - \widehat{H}_k}{\sqrt{\varepsilon}}\right) + o_P(1). \tag{28}$$

Using the fact that $\widehat{B}_k$ is a weakly consistent estimator of $b_k$ it is enough to prove the following convergence

$$\frac{(\sigma^k_{2H} - \sigma^k_{2\widehat{H}_k})}{\sqrt{\varepsilon}\log(\varepsilon)} \xrightarrow[\varepsilon \to 0]{\mathcal{P}} 0,$$

which is the same as showing

$$\frac{(\sigma^2_{2H} - \sigma^2_{2\widehat{H}_k})}{\sqrt{\varepsilon}\log(\varepsilon)} \xrightarrow[\varepsilon \to 0]{\mathcal{P}} 0. \tag{29}$$

We write

$$\frac{(\sigma_{2H}^2 - \sigma_{2\widehat{H}_k}^2)}{\sqrt{\varepsilon}\log(\varepsilon)} = \frac{1}{2\pi\sqrt{\varepsilon}\log(\varepsilon)} \int_{-\infty}^{\infty} |x|^{(3-2H)}|\widehat{\varphi}(-x)|^2 \{1 - \exp(2(H - \widehat{H}_k)\log(|x|))\}\, dx.$$

Making a Taylor expansion for the exponential function we obtain

$$\frac{(\sigma_{2H}^2 - \sigma_{2\widehat{H}_k}^2)}{\sqrt{\varepsilon}\log(\varepsilon)} = \frac{(\widehat{H}_k - H)}{\pi\sqrt{\varepsilon}\log(\varepsilon)} \int_{-\infty}^{\infty} |x|^{(3-2H)}|\widehat{\varphi}(-x)|^2 \log(|x|)\exp(\theta_\varepsilon(x))\, dx,$$

where $\theta_\varepsilon(x)$ is a point between 0 and $2(H - \widehat{H}_k)\log(|x|)$. By using assertion (1) of Theorem 3.2, and inequality

$$\exp(\theta_\varepsilon(x)) \le \exp(2|H - \widehat{H}_k||\log(|x|)|),$$

we will get the convergence in (29) by showing the following convergence result

$$\int_{-\infty}^{\infty} |x|^{(3-2H)}|\widehat{\varphi}(-x)|^2 |\log(|x|)|\exp(2|H - \widehat{H}_k||\log(|x|)|)\, dx$$

$$\xrightarrow[\varepsilon \to 0]{\text{a.s.}} \int_{-\infty}^{\infty} |x|^{(3-2H)}|\widehat{\varphi}(-x)|^2 |\log(|x|)|\, dx < +\infty. \tag{30}$$

To prove the convergence in (30) we use the fact that $(H - \widehat{H}_k) = o_{\text{a.s.}}(1)$ and then for $x \ne 0$,

$$|x|^{(3-2H)}|\widehat{\varphi}(-x)|^2 |\log(|x|)|\exp(2|H - \widehat{H}_k||\log(|x|)|)$$

$$\xrightarrow[\varepsilon \to 0]{\text{a.s.}} |x|^{(3-2H)}|\widehat{\varphi}(-x)|^2 |\log(|x|)|.$$

Now, let $0 < \delta < \inf(2H, (4 - 2H))$, then almost surely for all $\omega$, there exists $\varepsilon(\omega)$ such that $2|(H - \widehat{H}_k(\omega))| \le \delta$, when $\varepsilon \le \varepsilon(\omega)$. Furthermore, using the fact that $\varphi$ is a density one has $|\widehat{\varphi}(-x)|^2 \le 1$. Since for $x \ne 0$, $|\widehat{\varphi}(-x)|^2 \le \mathbf{C}|x|^{-4}$, then almost surely for all $\omega$, for all $\varepsilon \le \varepsilon(\omega)$ and $x \ne 0$, one obtains

$$|x|^{(3-2H)}|\widehat{\varphi}(-x)|^2 |\log(|x|)|\exp(2|H - \widehat{H}_k||\log(|x|)|)$$

$$\le |x|^{(3-2(H+\delta/2))}|\log(|x|)|\mathbf{1}_{|x|\le 1} + |x|^{(-1-2(H-\delta/2))}|\log(|x|)|\mathbf{1}_{|x|>1}.$$

Since $(H - \delta/2) > 0$ and $(4 - 2(H + \delta/2)) > 0$, we can apply Lebesgue's dominated convergence theorem to obtain the convergence in (30), thus we proved the convergence in (29) and equality (28) follows.

Now if we remark that the asymptotic behavior of $(\widehat{\sigma}_k/\sigma - 1)$ (see (17) for definition of $\widehat{\sigma}_k$) is the same as that of $(\widehat{\sigma^k}/\sigma^k - 1)/k$, then by (28) and assertion (1) of Theorem 3.2 the proof of assertion (2) is complete.

**Remark.** The last step of the proof shows that the asymptotic behavior of $((\hat{H}_k - H)/\sqrt{\varepsilon}, (\hat{\sigma}_k - \sigma)/(\sqrt{\varepsilon}\log(\varepsilon)))$ is a degenerated Gaussian law. □

**Proof of Theorem 3.3.**

(1) Assertion (1) follows from assertion (1) of Theorem 3.1.
(2) Assertion (2) of Theorem 3.1 and the remark following Theorem 2.2 imply that $\frac{1}{\sqrt{\varepsilon}}([\frac{\tilde{\sigma}_k}{\sigma}]^k - 1)$ converges weakly to $\sigma_{g_k}\mathcal{N}(0,1)$ that yields assertion (2).

Remark 2 follows from the fact that since $\hat{g}_{2,k} = \frac{k}{2}$, one has

$$\frac{\sigma_{g_k}^2}{k^2} = \frac{1}{k^2}\sum_{n=1}^{\infty} \hat{g}_{2n,k}^2 (2n)! \int_{-\infty}^{+\infty} \rho_H^{2n}(x)\, dx \ge \frac{2}{k^2}\hat{g}_{2,k}^2 \int_{-\infty}^{+\infty} \rho_H^2(x)\, dx = \frac{1}{2}\int_{-\infty}^{+\infty} \rho_H^2(x)\, dx = \frac{\sigma_{g_2}^2}{4}.$$

□

### 3.2. Functional estimation of $\sigma(\cdot)$

Under certain regularity conditions for $\mu(\cdot)$ and $\sigma(\cdot)$, we consider the "pseudo-diffusion" equations (1) with respect to $b_H(\cdot)$, that is

$$X(t) = c + \int_0^t \sigma(X(u))\,\mathrm{d}b_H(u) + \int_0^t \mu(X(u))\,\mathrm{d}u,$$

for $t \geq 0$, $H > 1/2$ and positive $\sigma(\cdot)$. We consider the problem of estimating $\sigma(\cdot)$. Suppose we observe, as before, instead of $X(t)$ the smoothed by convolution process $X_\varepsilon(t) = \frac{1}{\varepsilon}\int_{-\infty}^{\infty} \varphi(\frac{t-x}{\varepsilon})X(x)\,\mathrm{d}x$, with $\varphi(\cdot)$ as in Section 2, where we have extended $X(\cdot)$ by means of $X(t) = c$, if $t < 0$.

In a previous paper [2], in the case where $\mu(\cdot) \equiv 0$, estimation of $\sigma(\cdot)$ is done, using the first increments of $X(\cdot)$ or more generally the first derivative of $X_\varepsilon(\cdot)$. Namely we proved the two following theorems.

**Theorem 3.4.** *Let $1/2 < H < 1$. If $h(\cdot) \in C^0$ and $\sigma(\cdot) \in C^1$ then*

$$\sqrt{\frac{\pi}{2}}\frac{\varepsilon^{(1-H)}}{\tilde{\sigma}_{2H}} \int_0^1 h(X_\varepsilon(u))|\dot{X}_\varepsilon(u)|\,\mathrm{d}u \xrightarrow{a.s.} \int_0^1 h(X(u))\sigma(X(u))\,\mathrm{d}u,$$

*where $\tilde{\sigma}_{2H}$ is defined by*

$$\tilde{\sigma}_{2H}^2 := \mathbb{V}[\varepsilon^{(1-H)}\dot{b}_H^\varepsilon(t)] = \frac{1}{2\pi}\int_{-\infty}^{+\infty} |x|^{(1-2H)}|\hat{\varphi}(-x)|^2\,\mathrm{d}x.$$

The rate of convergence is given by Theorem 3.5.

**Theorem 3.5.** *Let us suppose that $1/2 < H < 3/4$, $h(\cdot) \in C^4$, $\sigma(\cdot) \in C^4$, $\sigma(\cdot)$ is bounded and $\sup\{|\sigma^{(4)}(x)|, |h^{(4)}(x)|\} \leq P(|x|)$, where $P(\cdot)$ is a polynomial, then*

$$\frac{1}{\sqrt{\varepsilon}}\left[\sqrt{\frac{\pi}{2}}\frac{\varepsilon^{(1-H)}}{\tilde{\sigma}_{2H}} \int_0^1 h(X_\varepsilon(u))|\dot{X}_\varepsilon(u)|\,\mathrm{d}u - \int_0^1 h(X(u))\sigma(X(u))\,\mathrm{d}u\right],$$

*converges stably towards*

$$\sigma_{g_1}\int_0^1 h(X(u))\sigma(X(u))\,\mathrm{d}\widehat{W}(u).$$

*Here, $\widehat{W}(\cdot)$ is a standard Brownian motion independent of $b_H(\cdot)$, $g_1(x) = \sqrt{\frac{\pi}{2}}|x| - 1$.*

We give an outline of the proof of the last two above theorems in order to generalize these results to our setting, considering the case where $\mu(\cdot)$ is not necessarily null. Indeed, because $b_H(\cdot)$ has zero quadratic variation when $H > 1/2$, Lin (see ([8])) proved that when $\sigma(\cdot) \in C^1$ and $\mu(\cdot) \equiv 0$, the solution to the stochastic differential equation (1) can be expressed as $X(t) = K(b_H(t))$, for $t \geq 0$, where $K(t)$ is the solution to the ordinary differential equation (ODE)

$$\dot{K}(t) = \sigma(K(t)), \qquad K(0) = c \tag{31}$$

(for $t < 0$, $X(t) = c$).

**Heuristic proof of Theorem 3.4.** We have shown in [2] that

$$\sqrt{\frac{\pi}{2}}\frac{\varepsilon^{(1-H)}}{\tilde{\sigma}_{2H}} \int_0^1 h(X_\varepsilon(u))|\dot{X}_\varepsilon(u)|\,\mathrm{d}u$$

$$\simeq \sqrt{\frac{\pi}{2}}\frac{\varepsilon^{(1-H)}}{\tilde{\sigma}_{2H}} \int_0^1 h(K(b_H^\varepsilon(u)))\sigma(K(b_H^\varepsilon(u)))|\dot{b}_H^\varepsilon(u)|\,\mathrm{d}u,$$

and then Theorem 3.4 ensues from the Azaïs and Wschebor Theorem (see [1]) that follows.

**Theorem 3.6.** *Let $1/2 < H < 1$. For every continuous function $h(\cdot)$*

$$\sqrt{\frac{\pi}{2}} \frac{\varepsilon^{(1-H)}}{\tilde{\sigma}_{2H}} \int_0^1 h(b_H^\varepsilon(u))|\dot{b}_H^\varepsilon(u)|\,\mathrm{d}u \xrightarrow{a.s.} \int_0^1 h(b_H(u))\,\mathrm{d}u.$$

$\square$

**Heuristic proof of Theorem 3.5.** We proved in [2] that the cited approximation is $o_{\mathrm{a.s.}}(\sqrt{\varepsilon})$ that is

$$\varepsilon^{(1-H)}\left[\int_0^1 h(X_\varepsilon(u))|\dot{X}_\varepsilon(u)|\,\mathrm{d}u - \int_0^1 h(K(b_H^\varepsilon(u)))\sigma(K(b_H^\varepsilon(u)))|\dot{b}_H^\varepsilon(u)|\,\mathrm{d}u\right] = o_{\mathrm{a.s.}}(\sqrt{\varepsilon}),$$

hence the asymptotic behavior of

$$\frac{1}{\sqrt{\varepsilon}}\left[\sqrt{\frac{\pi}{2}}\frac{\varepsilon^{(1-H)}}{\tilde{\sigma}_{2H}}\int_0^1 h(X_\varepsilon(u))|\dot{X}_\varepsilon(u)|\,\mathrm{d}u - \int_0^1 h(X(u))\sigma(X(u))\,\mathrm{d}u\right],$$

is the same as that of

$$\frac{1}{\sqrt{\varepsilon}}\left[\frac{\varepsilon^{(1-H)}}{\tilde{\sigma}_{2H}}\sqrt{\frac{\pi}{2}}\int_0^1 h(K(b_H^\varepsilon(u)))\sigma(K(b_H^\varepsilon(u)))|\dot{b}_H^\varepsilon(u)|\,\mathrm{d}u - \int_0^1 h(K(b_H(u)))\sigma(K(b_H(u)))\,\mathrm{d}u\right],$$

hence Theorem 3.5 arises from Theorem 3.7 (see [2]) that follows.

**Theorem 3.7.** *Let $h(\cdot) \in C^4$ such that $|h^{(4)}(x)| \leq P(|x|)$, where $P(\cdot)$ is a polynomial, if $1/4 < H < 3/4$, then*

$$\frac{1}{\sqrt{\varepsilon}}\left[\sqrt{\frac{\pi}{2}}\frac{\varepsilon^{(1-H)}}{\tilde{\sigma}_{2H}}\int_0^1 h(b_H^\varepsilon(u))|\dot{b}_H^\varepsilon(u)|\,\mathrm{d}u - \int_0^1 h(b_H(u))\,\mathrm{d}u\right],$$

*converges stably towards*

$$\sigma_{g_1}\int_0^1 h(b_H(u))\,\mathrm{d}\widehat{W}(u).$$

*Here, $\widehat{W}(\cdot)$ is still a standard Brownian motion independent of $b_H(\cdot)$, $g_1(x) = \sqrt{\frac{\pi}{2}}|x| - 1$.*

$\square$

Now we get back to our purpose, i.e. assuming that $\mu(\cdot)$ is not necessarily null in model (1), and estimating $\sigma(\cdot)$ by considering the second-order increments of $X(\cdot)$ and more generally the second derivative of $X_\varepsilon(\cdot)$ and working with $|\ddot{X}_\varepsilon(\cdot)|^k$ with $k \geq 1$, instead of $|\dot{X}_\varepsilon(\cdot)|$.

In this aim, we are considering the following assumptions on the coefficients $\mu(\cdot)$ and $\sigma(\cdot)$:

(H1) $\sigma(\cdot)$ is $C^1$ and Lipschitz function on $\mathbb{R}$, and bounded away from zero.

There exists some constant $\eta$, $(1/H - 1) < \eta \leq 1$, and for every $N > 0$, there exists $M_N > 0$ such that

$$|\dot{\sigma}(x) - \dot{\sigma}(y)| \leq M_N|x - y|^\eta, \quad \forall |x|, |y| \leq N.$$

(H2) $\mu(\cdot)$ is $C^1$, bounded and Lipschitz function on $\mathbb{R}$.

***Remark.*** *Hypotheses (H1) and (H2) require that $\sigma(\cdot)$ is bounded away from zero and that $\mu(\cdot)$ is bounded. These last two assumptions can be replaced by the following: there exist $0 \leq \gamma \leq 1$ and $M > 0$ such that, $|\sigma(x)| \leq M(1 + |x|^\gamma)$ and $|\mu(x)| \leq M(1 + |x|)$, for all $x \in \mathbb{R}$. These new assumptions ensure that there exists an unique process solution to the stochastic equation (1). Furthermore $X(\cdot)$ has almost surely $(H - \delta)$-Hölder continuous trajectories on all compact included in $\mathbb{R}^+$ (see [11]). In particular, the solutions $X(\cdot)$ to the last four previous models in Section 3.1 are almost surely $(H - \delta)$-Hölder continuous.*

At this step, let us state two theorems that we get thanks to Theorem 3.6 and Theorem 3.7.

**Theorem 3.8.** *If $h(\cdot) \in C^0$ and $1/2 < H < 1$, under hypotheses* (H1) *and* (H2) *then,*

$$\frac{1}{\mathbb{E}[|N|^k]} \int_0^1 h(X_\varepsilon(u)) \left| \frac{\varepsilon^{(2-H)} \ddot{X}_\varepsilon(u)}{\sigma_{2H}} \right|^k du \xrightarrow{a.s.} \int_0^1 h(X(u)) [\sigma(X(u))]^k du.$$

***Remark.*** *If $\mu(\cdot) \equiv 0$, hypotheses* (H1) *and* (H2) *can be replaced by $\sigma(\cdot) \in C^1$.*

**Theorem 3.9.** *Let us suppose that $1/2 < H < 1$, $h(\cdot) \in C^4$, $\sigma(\cdot) \in C^4$, $\sigma(\cdot)$ is bounded and $\sup\{|\sigma^{(4)}(x)|, |h^{(4)}(x)|\} \leq P(|x|)$, where $P(\cdot)$ is a polynomial, then under hypotheses* (H1) *and* (H2)

$$\frac{1}{\sqrt{\varepsilon}} \left[ \frac{1}{\mathbb{E}[|N|^k]} \int_0^1 h(X_\varepsilon(u)) \left| \frac{\varepsilon^{(2-H)}}{\sigma_{2H}} \ddot{X}_\varepsilon(u) \right|^k du - \int_0^1 h(X(u)) [\sigma(X(u))]^k du \right]$$

$$\underset{\varepsilon \to 0}{\Longrightarrow} \sigma_{g_k} \int_0^1 h(X(u)) [\sigma(X(u))]^k d\widehat{W}(u).$$

*Here, $g_k(x) = \frac{1}{\mathbb{E}[|N|^k]} |x|^k - 1$ and $\widehat{W}(\cdot)$ is still a standard Brownian motion independent of $b_H(\cdot)$.*

***Remark 1.*** *If $\mu(\cdot) \equiv 0$, hypotheses* (H1) *and* (H2) *can be relaxed and the convergence will be stably convergence.*

***Remark 2.*** *Note that this theorem is valid for all the range of $H$ belonging to $]1/2, 1[$ instead of $]1/2, 3/4[$ as in Theorem 3.5, this is due to the fact that we are working with the second-order increments of $X(\cdot)$ instead of the first ones.*

***Remark 3.*** *Note that by Theorem 3.1 and by the remark following Theorem 2.2, Theorem 3.8 and Theorem 3.9 are still available under hypothesis $h(\cdot) \equiv 1$ for the first two models, even if hypotheses* (H1) *and* (H2) *are not exactly satisfied by the second model.*

Although these hypotheses are not fulfilled by the third model these two last results of convergence remain valid, i.e. we have that

$$\frac{1}{\mathbb{E}[|N|^k]} \int_0^1 \left| \frac{\varepsilon^{(2-H)} \ddot{X}_\varepsilon(u)}{\sigma_{2H}} \right|^k du \xrightarrow{a.s.} \sigma^k \int_0^1 |X(u)|^k du$$

and that

$$\frac{1}{\sqrt{\varepsilon}} \left[ \frac{1}{\mathbb{E}[|N|^k]} \int_0^1 \left| \frac{\varepsilon^{(2-H)}}{\sigma_{2H}} \ddot{X}_\varepsilon(u) \right|^k du - \sigma^k \int_0^1 |X(u)|^k du \right]$$

converges weakly towards $\sigma_{g_k} \sigma^k \int_0^1 |X(u)|^k d\widehat{W}(u)$.

***Remark 4.*** *We conjecture that the two last results of convergence are still true for the fourth model, even if hypotheses* (H1) *and* (H2) *are not satisfied by this one.*

**Proof of Theorems 3.8 and 3.9.** We begin by showing the remark (resp. Remark 1) following Theorem 3.8 (resp. Theorem 3.9) i.e. we suppose $\mu(\cdot) \equiv 0$ in model (1).

We know that since $\sigma(\cdot) \in C^1$, for $t \geq 0$, $X(t) = K(b_H(t))$, where $K(\cdot)$ is the solution to the ODE (31) and for $t < 0$, $X(t) = c$. Thus we are going to prove that

$$T_\varepsilon(h) := \left[ \int_0^1 h(X_\varepsilon(u)) |\varepsilon^{(2-H)} \ddot{X}_\varepsilon(u)|^k du - \int_0^1 h(K(b_H^\varepsilon(u))) [\sigma(K(b_H^\varepsilon(u)))]^k |\varepsilon^{(2-H)} \ddot{b}_H^\varepsilon(u)|^k du \right]$$

is $o_{a.s.}(\sqrt{\varepsilon})$ when $h(\cdot) \in C^1$ and $o_{a.s.}(1)$ when $h(\cdot) \in C^0$. Then the remarks will follow from a generalization of Theorems 3.6 and 3.7. The proofs of these generalizations being easy to obtain, they will not be given here. We still just remark that since instead of considering the first-order increments of $b_H(\cdot)$, we study the second-order ones, and convergence in the generalization of Theorem 3.7 is reached for all values of $H$ in $]1/2, 1[$.

We need the following lemma for which a proof is provided in the Appendix.

**Lemma 3.2.** *In model* (1), *for* $\mu(\cdot) \equiv 0$ *and for* $0 \leq t \leq 1$,

$$X_\varepsilon(t) = K(b_H^\varepsilon(t)) + a_\varepsilon(t) \quad and \quad \ddot{X}_\varepsilon(t) = \sigma(K(b_H^\varepsilon(t)))\ddot{b}_H^\varepsilon(t) + c_\varepsilon(t),$$

*with*

$$|a_\varepsilon(t)| \leq \mathbf{C}(\omega)\varepsilon^{(H-\delta)},$$

*and*

$$|c_\varepsilon(t)| \leq \mathbf{C}(\omega)(\varepsilon^{(H-2-\delta)}\mathbf{1}_{\{0 \leq t \leq \varepsilon\}} + \varepsilon^{(2H-2-\delta)}\mathbf{1}_{\{\varepsilon \leq t \leq 1\}}) \quad \text{for any } \delta > 0.$$

Now,

$$T_\varepsilon(h) = L_1 + L_2,$$

where

$$L_1 := \int_0^1 [h(X_\varepsilon(u)) - h(K(b_H^\varepsilon(u)))]|\varepsilon^{(2-H)}\ddot{X}_\varepsilon(u)|^k \, du$$

and

$$L_2 := \int_0^1 h(K(b_H^\varepsilon(u)))\{|\varepsilon^{(2-H)}\ddot{X}_\varepsilon(u)|^k - |\sigma(K(b_H^\varepsilon(u)))\varepsilon^{(2-H)}\ddot{b}_H^\varepsilon(u)|^k\} \, du.$$

Now, let us study $L_1$ and $L_2$. For $L_1$, if $h(\cdot) \in C^1$, we have

$$L_1 = \int_0^1 \dot{h}(\theta)(X_\varepsilon(u) - K(b_H^\varepsilon(u)))|\varepsilon^{(2-H)}\ddot{X}_\varepsilon(u)|^k \, du,$$

where $\theta$ is a point between $X_\varepsilon(u)$ and $K(b_H^\varepsilon(u))$ and then by Lemma 3.2

$$|L_1| \leq \mathbf{C}(\omega) \int_0^1 |X_\varepsilon(u) - K(b_H^\varepsilon(u))||\varepsilon^{(2-H)}\ddot{X}_\varepsilon(u)|^k \, du$$

$$\leq \mathbf{C}(\omega)\varepsilon^{(H-\delta)} \int_0^1 |\varepsilon^{(2-H)}\ddot{X}_\varepsilon(u)|^k \, du = o_{a.s.}(\sqrt{\varepsilon}),$$

since $H > 1/2$ and because of the boundness of $\int_0^1 |\varepsilon^{(2-H)}\ddot{X}_\varepsilon(u)|^k \, du$.

This last remark can be shown by writing

$$\|\varepsilon^{(2-H)}\ddot{X}_\varepsilon(\cdot)\|_k \leq \varepsilon^{(2-H)}\|\ddot{X}_\varepsilon(\cdot) - \sigma(K(b_H^\varepsilon(\cdot)))\ddot{b}_H^\varepsilon(\cdot)\|_k + \|\sigma(K(b_H^\varepsilon(\cdot)))\varepsilon^{(2-H)}\ddot{b}_H^\varepsilon(\cdot)\|_k,$$

and then by Lemma 3.2 and since $\int_0^1 |\varepsilon^{(2-H)}\ddot{b}_H^\varepsilon(u)|^k \, du$ is bounded (see the remark following Theorem 2.1), one gets for $\delta$ small enough

$$\|\varepsilon^{(2-H)}\ddot{X}_\varepsilon(\cdot)\|_k \leq \mathbf{C}(\omega)(\varepsilon^{(1/k-\delta)} + \varepsilon^{(H-\delta)} + 1) \leq \mathbf{C}(\omega).$$

Note that if $h(\cdot)$ is only in $C^0$, the last remark and Lemma 3.2 imply that $L_1$ tends almost surely to zero when $\varepsilon$ goes to zero.

Moreover, for $L_2$ we have

$$|L_2| \leq \mathbf{C}(\omega) \| |\varepsilon^{(2-H)} \ddot{X}_\varepsilon(\cdot)|^k - |\sigma(K(b_H^\varepsilon(\cdot))) \varepsilon^{(2-H)} \ddot{b}_H^\varepsilon(\cdot)|^k \|_1. \tag{32}$$

Furthermore, applying Lemma 3.2 we obtain

$$\varepsilon^{(2-H)} \| \ddot{X}_\varepsilon(\cdot) - \sigma(K(b_H^\varepsilon(\cdot))) \ddot{b}_H^\varepsilon(\cdot) \|_{k,\varepsilon} \leq \mathbf{C}(\omega) \varepsilon^{(1/k-\delta)}$$

and

$$\varepsilon^{(2-H)} \| \ddot{X}_\varepsilon(\cdot) - \sigma(K(b_H^\varepsilon(\cdot))) \ddot{b}_H^\varepsilon(\cdot) \|_{k,\varepsilon^c} \leq \mathbf{C}(\omega) \varepsilon^{(H-\delta)}.$$

Now as in the proof of Theorem 3.1 we get $\|\varepsilon^{(2-H)} \ddot{b}_H^\varepsilon(\cdot)\|_{k,\varepsilon} \leq \mathbf{C}(\omega) \varepsilon^{(1/k-\delta)}$ and $\|\varepsilon^{(2-H)} \ddot{b}_H^\varepsilon(\cdot)\|_{k,\varepsilon^c} \leq \mathbf{C}(\omega) \varepsilon^{-\delta}$.

Finally taking $\delta$ small enough, i.e. $0 < \delta < (H - 1/2)/k$, we prove applying the second part of inequality (18) to inequality (32) that

$$|L_2| \leq \mathbf{C}(\omega)(\varepsilon^{(1-\delta k)} + \varepsilon^{(H-\delta k)} + \varepsilon^{k(H-\delta)}) \leq \mathbf{C}(\omega) \varepsilon^{(H-\delta k)} = \mathrm{o}_{\mathrm{a.s.}}(\sqrt{\varepsilon}),$$

and this inequality completes the proof, i.e. we have shown that $T_\varepsilon(h) = \mathrm{o}_{\mathrm{a.s.}}(\sqrt{\varepsilon})$ when $h(\cdot) \in C^1$ and $\mathrm{o}_{\mathrm{a.s.}}(1)$ when $h(\cdot) \in C^0$.

Now we get back to the model (1) with $\mu(\cdot)$ not necessarily identically null and we are going to prove Theorem 3.8 (resp. Theorem 3.9) using the remark (resp. Remark 1) following it.

Let $X(t)$ be the solution to the equation

$$\mathrm{d}X(t) = \sigma(X(t)) \, \mathrm{d}b_H(t) + \mu(X(t)) \, \mathrm{d}t.$$

We denote $P$ the probability measure induced by the fBm over the $\sigma$-algebra $\mathcal{G}$. If $G$ is a measurable and bounded real function defined on the space $C([0,1], \mathbb{R})$ of continuous real functions, we have

$$\mathbb{E}[G(X)]_P = \mathbb{E}[G(K(b_H)) \Lambda]_P, \tag{33}$$

where $\Lambda$ will be defined later on and $K(\cdot)$ is the solution to the ODE (31). To obtain this equality, we use hypotheses (H1) and (H2) to apply the Girsanov theorem of Decreusefond–Üstünel (see [6]). Namely let $Y(t) = K(b_H(t))$ and let define $\tilde{b}_H(t) := b_H(t) - \int_0^t \frac{\mu(Y(s))}{\sigma(Y(s))} \, \mathrm{d}s$. By using the Itô's formula we get

$$\mathrm{d}Y(t) = \sigma(Y(t)) \, \mathrm{d}\tilde{b}_H(t) + \mu(Y(t)) \, \mathrm{d}t.$$

Furthermore, there exists a probability measure $\tilde{P}$ absolutely continuous w.r.t. $P$, such that with this probability $\tilde{P}$ the process $\tilde{b}_H(\cdot)$ is a fBm with parameter $0 < H < 1$. Hence, we have

$$\mathbb{E}[G(Y)]_{\tilde{P}} = \mathbb{E}[G(K(b_H)) \Lambda]_P,$$

where $\Lambda$ is the Radon-Nikodym derivative of $\tilde{P}$ w.r.t. $P$. Since the two processes $X(\cdot)$ and $Y(\cdot)$ have the same distribution over $P$ and $\tilde{P}$ respectively, we get

$$\mathbb{E}[G(X)]_P = \mathbb{E}[G(Y)]_{\tilde{P}} = \mathbb{E}[G(K(b_H)) \Lambda]_P,$$

and equality (33) follows.

Let us define the set of trajectories

$$\Delta := \left\{ x \in C([0,1], \mathbb{R}) : \lim_{\varepsilon \to 0} \frac{1}{\mathbb{E}[|N|^k]} \int_0^1 h(x_\varepsilon(u)) \left| \frac{\varepsilon^{(2-H)} \ddot{x}_\varepsilon(u)}{\sigma_{2H}} \right|^k \mathrm{d}u = \int_0^1 h(x(u)) [\sigma(x(u))]^k \, \mathrm{d}u \right\}.$$

If we choose $G$ as $1_\Delta$, using (33) one obtains

$$\mathbb{E}[1_\Delta(X)]_P = \mathbb{E}[1_\Delta(K(b_H))\Lambda]_P = \mathbb{E}[\Lambda]_P = 1,$$

where above we have used the remark following Theorem 3.8 i.e. $P(K(b_H) \in \Delta) = 1$, thus Theorem 3.8 follows.

In what follows we want to show Theorem 3.9 and in this aim we must study the weak convergence in Theorem 3.8. Consider

$$M_\varepsilon(Y) = \frac{1}{\sqrt{\varepsilon}} \left[ \frac{1}{\mathbb{E}[|N|^k]} \int_0^1 h(Y_\varepsilon(u)) \left| \frac{\varepsilon^{(2-H)}}{\sigma_{2H}} \ddot{Y}_\varepsilon(u) \right|^k du - \int_0^1 h(Y(u))[\sigma(Y(u))]^k du \right],$$

by Remark 1 following Theorem 3.9 we know that this term stably converges towards

$$M(Y, \widehat{W}) = \sigma_{g_k} \int_0^1 h(Y(u))[\sigma(Y(u))]^k d\widehat{W}(u).$$

Once again we use the Girsanov theorem. Namely, let $F$ be a continuous and bounded real function, then by applying equality (33) to $G = F \circ M_\varepsilon$, one gets

$$\mathbb{E}[F(M_\varepsilon(X))]_P = \mathbb{E}[F(M_\varepsilon(K(b_H)))\Lambda]_P,$$

given that $\Lambda$ is a measurable function w.r.t. the $\sigma$-algebra $\mathcal{G}$, and using the properties of stable convergence we get

$$\mathbb{E}[F(M_\varepsilon(K(b_H)))\Lambda]_P \to \mathbb{E}[F(M(Y, \widehat{W}))\Lambda]_{P \otimes P_1},$$

where $\mathbb{E}[\cdot]_{P \otimes P_1}$ denotes the expectation w.r.t. the product probability of the fBm and the independent Brownian motion $\widehat{W}(\cdot)$. By using the Girsanov theorem another time we have

$$E[F(M(X, \widehat{W}))]_{P \otimes P_1} = E[F(M(Y, \widehat{W}))\Lambda]_{P \otimes P_1}.$$

To prove the last equality, suppose that $M(Y, \widehat{W})$ depends only on a finite number of coordinates with respect to the second variable, namely:

$$M(Y, \widehat{W}) = M(Y, \widehat{W}(t_1), \widehat{W}(t_2), \ldots, \widehat{W}(t_m)),$$

thus if we note $p_{t_1,t_2,\ldots,t_m}(x_1, x_2, \ldots, x_m)$ the density of the vector $(\widehat{W}(t_1), \widehat{W}(t_2), \ldots, \widehat{W}(t_m))$, we have by independence

$$\mathbb{E}[F(M(Y, \widehat{W}(t_1), \widehat{W}(t_2), \ldots, \widehat{W}(t_m)))\Lambda]_{P \otimes P_1}$$
$$= \int_{R^m} \mathbb{E}[F(M(Y, x_1, x_2, \ldots, x_m))\Lambda]_P p_{t_1,t_2,\ldots,t_m}(x_1, x_2, \ldots, x_m) dx_1 dx_2 \cdots dx_m. \tag{34}$$

Applying formula (33) to $G = F \circ M(\cdot, x_1, x_2, \ldots, x_m)$ in equality (34) we get

$$\mathbb{E}[F(M(Y, \widehat{W}(t_1), \widehat{W}(t_2), \ldots, \widehat{W}(t_m)))\Lambda]_{P \otimes P_1}$$
$$= \int_{R^m} \mathbb{E}[F(M(X, x_1, x_2, \ldots, x_m))]_P p_{t_1,t_2,\ldots,t_m}(x_1, x_2, \ldots, x_m) dx_1 dx_2 \cdots dx_m$$
$$= \mathbb{E}[F(M(X, \widehat{W}(t_1), \widehat{W}(t_2), \ldots, \widehat{W}(t_m)))]_{P \otimes P_1},$$

by using independence again.

A classical approximation argument leads to the conclusion. Hence we proved that

$$\mathbb{E}[F(M_\varepsilon(X))]_P \to \mathbb{E}[F(M(X,\widehat{W}))]_{P\otimes P_1},$$

and Theorem 3.9 follows.

Remark 3 requires the following argumentation.

Hypotheses (H1) and (H2) on $\mu(\cdot)$ and $\sigma(\cdot)$ are given here to apply the Girsanov theorem but they are somewhat restrictive. Thus, the third model does not verify these hypotheses, nevertheless, the Girsanov theorem can be applied to this model. Indeed, as before let define $\tilde{b}_H(t) := b_H(t) - \int_0^t \frac{\mu(Y(s))}{\sigma(Y(s))}\,\mathrm{d}s = b_H(t) - \frac{\mu}{\sigma}t$. There exists $\tilde{P}$, a probability measure absolutely continuous w.r.t. $P$ such that over this probability the process $\tilde{b}_H(\cdot)$ is a fBm with parameter $0 < H < 1$.

To show this last statement, it is sufficient to prove that

$$\mathbb{E}\left[\exp\left(\frac{1}{2}\|\xi\|_{\mathcal{H}_H}^2\right)\right] < +\infty,$$

with $\xi_t = \int_0^t \frac{\mu(Y(s))}{\sigma(Y(s))}\,\mathrm{d}s = \frac{\mu}{\sigma}t$ (see Lemma 6 of [6] and Theorem 4.9 of [5] for notations and details of the proof of this argument).

To prove the finiteness of last expectation we use, as in Lemma 6 of [6], the following upper bound $\|\xi\|_{\mathcal{H}_H}^2 \leq \mathbf{C}\|\xi\|_{B_{2,2}^{H+1/2}}^2$. This last norm is equivalent to the Sobolev norm hence $\|\xi\|_{B_{2,2}^{H+1/2}}^2 \leq \mathbf{C}(\|\xi\|_2^2 + \|\dot{\xi}\|_{B_{2,2}^{H-1/2}}^2)$, and the two terms in the last right-hand expression are bounded independently of $b_H(\cdot)$.

To get Theorems 3.8 and 3.9 for the third model and for $h(\cdot) \equiv 1$, it remains to verify that these theorems are valid in the case where $h(\cdot) \equiv 1$ and where $X(\cdot)$ is the solution of (6) with $\mu(\cdot) = 0$. For this we need to apply generalizations of Theorem 3.6 and of Theorem 3.7 in the case where $h(x) = |\sigma(K(x))|^k = \sigma^k|c|^k\exp(k\sigma x)$. Thus, since $h(\cdot)$ is a continuous function Theorem 3.8 follows. Nevertheless even if $h^{(4)}(\cdot)$ cannot be bounded by a polynomial, it can be shown that Theorem 3.7 is still true. Finally Remark 3 is fulfilled. □

### 3.3. Tests of the hypothesis

#### 3.3.1. Three simple models

Let us consider the three stochastic differential equations, for $t \geq 0$,

$$\mathrm{d}X_\varepsilon(t) = \sigma_\varepsilon(X_\varepsilon(t))\,\mathrm{d}b_H(t) + \mu(X_\varepsilon(t))\,\mathrm{d}t, \quad \text{with } X_\varepsilon(0) = c, \tag{35}$$

$X_\varepsilon(t) = c$, for $t < 0$ and where $H$, $\frac{1}{2} < H < 1$, is known. We consider testing the hypothesis

$$H_0: \quad \sigma_\varepsilon(x) = \sigma \quad (\text{resp. } \sigma_\varepsilon(x) = \sigma x),$$

against the sequence of alternatives

$$H_\varepsilon: \quad \sigma_\varepsilon(x) = \sigma_\varepsilon := \sigma + \sqrt{\varepsilon}(d + F(\sqrt{\varepsilon})) \quad (\text{resp. } \sigma_\varepsilon(x) = \sigma_\varepsilon x),$$

where $\sigma$, $d$ are positive constants, $F(\cdot)$ is a positive function such that $F(\sqrt{\varepsilon}) \xrightarrow[\varepsilon \to 0]{} 0$ and $\mu(x) = \mu$ or $\mu(x) = \mu x$ for the first two models (resp. $\mu(x) = \mu x$ for the third one).

By Section 3.1, conditions on $\sigma_\varepsilon(\cdot)$ and on $\mu(\cdot)$ ensure that for each model there exists an unique process solution to the stochastic equation (35), let us say $X_\varepsilon(\cdot)$. Let us suppose that the observed process is

$$Y_\varepsilon(\cdot) := \frac{1}{\varepsilon}\int_{-\infty}^{+\infty}\varphi\left(\frac{\cdot - x}{\varepsilon}\right)X_\varepsilon(x)\,\mathrm{d}x,$$

with $\varphi(\cdot)$ as in Section 2. We are interested in observing the following functionals

$$F_\varepsilon := \frac{1}{\sqrt{\varepsilon}}\left[\sqrt{\frac{\pi}{2}}\frac{\varepsilon^{(2-H)}}{\sigma_{2H}}\int_0^1 |\ddot{Y}_\varepsilon(u)|\,\mathrm{d}u - \sigma\right]$$

$$\left(\text{resp. } F_\varepsilon := \frac{1}{\sqrt{\varepsilon}}\left[\sqrt{\frac{\pi}{2}}\frac{\varepsilon^{(2-H)}}{\sigma_{2H}}\int_0^1 |\ddot{Y}_\varepsilon(u)|\,\mathrm{d}u - \sigma\int_0^1 |Y_\varepsilon(u)|\,\mathrm{d}u\right]\right).$$

Using generalizations of Theorems 3.6 and 3.7 we can prove the following theorem.

**Theorem 3.10.** *Let us suppose that $H$ is known with $1/2 < H < 1$, then*

$$F_\varepsilon \underset{\varepsilon\to 0}{\Longrightarrow} \sigma_{g_1}\sigma N + d$$

*(resp. $F_\varepsilon$ converges stably towards*

$$\sigma_{g_1}\sigma\int_0^1 |X(u)|\,\mathrm{d}\widehat{W}(u) + d\int_0^1 |X(u)|\,\mathrm{d}u,$$

*where $X(\cdot)$ is the solution to (6) and $\widehat{W}(\cdot)$ is a standard Brownian motion independent of $X(\cdot)$) and $g_1(\cdot)$ is defined by (10).*

**Remark 1.** There is an asymptotic bias $d$ (resp. a random asymptotic bias $d\int_0^1 |X(u)|\,\mathrm{d}u$), and the larger the bias the easier it is to discriminate between the two hypotheses.

**Remark 2.** $X_\varepsilon(\cdot)$ plays the role of $X(\cdot)$, $Y_\varepsilon(\cdot)$ that of $X_\varepsilon(\cdot)$ in first part of Section 3.1, with $\sigma_\varepsilon = \sigma$.

**Proof of Theorem 3.10.** We need the following lemma for which a proof is provided in the Appendix.

**Lemma 3.3.** *For $0 \le t \le 1$,*

$$\ddot{Y}_\varepsilon(t) = \sigma_\varepsilon \ddot{b}^\varepsilon_H(t) + a_\varepsilon(t)$$

$$(\text{resp. } \ddot{Y}_\varepsilon(t) = \sigma_\varepsilon X_\varepsilon(t)\ddot{b}^\varepsilon_H(t) + a_\varepsilon(t), \quad \text{and} \quad Y_\varepsilon(t) = X_\varepsilon(t) + d_\varepsilon(t)),$$

*where*

$$|a_\varepsilon(t)| \le \mathbf{C}(\omega)(\varepsilon^{(H-2-\delta)}\mathbf{1}_{\{0\le t\le \varepsilon\}} + \varepsilon^{(2H-2-\delta)}\mathbf{1}_{\{\varepsilon\le t\le 1\}})$$

*and*

$$|d_\varepsilon(t)| \le \mathbf{C}(\omega)\varepsilon^{(H-\delta)} \quad \text{for any } \delta > 0.$$

Now we write $F_\varepsilon$ as

$$F_\varepsilon = \sigma S_{g_1}(\varepsilon) + d + G_\varepsilon$$

$$\left(\text{resp. } F_\varepsilon = \frac{\sigma}{\sqrt{\varepsilon}}\int_0^1 |X(u)|g_1(Z_\varepsilon(u))\,\mathrm{d}u + d\int_0^1 |X(u)|\,\mathrm{d}u + G_\varepsilon\right),$$

where

$$G_\varepsilon := d\int_0^1 g_1(Z_\varepsilon(u))\,\mathrm{d}u + F(\sqrt{\varepsilon})\int_0^1 \sqrt{\frac{\pi}{2}}|Z_\varepsilon(u)|\,\mathrm{d}u + \frac{1}{\sqrt{\varepsilon}}\sqrt{\frac{\pi}{2}}\frac{\varepsilon^{(2-H)}}{\sigma_{2H}}\int_0^1 (|\ddot{Y}_\varepsilon(u)| - |\sigma_\varepsilon \ddot{b}^\varepsilon_H(u)|)\,\mathrm{d}u$$

$$\left( \text{resp. } G_\varepsilon := d \int_0^1 |X(u)| g_1(Z_\varepsilon(u)) \, \mathrm{d}u + F(\sqrt{\varepsilon}) \int_0^1 |X_\varepsilon(u)| \sqrt{\frac{\pi}{2}} |Z_\varepsilon(u)| \, \mathrm{d}u \right.$$

$$+ \frac{1}{\sqrt{\varepsilon}} \sqrt{\frac{\pi}{2}} \frac{\varepsilon^{(2-H)}}{\sigma_{2H}} \int_0^1 (|\ddot{Y}_\varepsilon(u)| - |\sigma_\varepsilon X_\varepsilon(u) \ddot{b}_H^\varepsilon(u)|) \, \mathrm{d}u$$

$$+ d \int_0^1 (|X_\varepsilon(u)| - |X(u)|) \sqrt{\frac{\pi}{2}} |Z_\varepsilon(u)| \, \mathrm{d}u + \frac{\sigma}{\sqrt{\varepsilon}} \int_0^1 (|X_\varepsilon(u)| - |Y_\varepsilon(u)|) \, \mathrm{d}u$$

$$\left. + \frac{\sigma}{\sqrt{\varepsilon}} \int_0^1 (|X_\varepsilon(u)| - |X(u)|) g_1(Z_\varepsilon(u)) \, \mathrm{d}u \right).$$

The remark following Theorem 2.1, Lemma 3.3 and the fact that $F(\sqrt{\varepsilon})$ tends to zero with $\varepsilon$ ensure that $G_\varepsilon = \mathrm{o}_{\mathrm{a.s.}}(1)$ (resp. the same arguments, the fact that $(|X_\varepsilon(u)| - |X(u)|)/\sqrt{\varepsilon}$ almost surely uniformly converges towards

$$|c| d \exp(\mu u + \sigma b_H(u)) b_H(u)$$

and a generalization of Theorem 3.6 give that $G_\varepsilon = \mathrm{o}_{\mathrm{a.s.}}(1)$).

The remark following Theorem 2.2 allows us to conclude (resp. the asymptotic behavior of $F_\varepsilon$ can be treated in the same manner that we have done in Theorem 3.7, where instead of working with $h(b_H(u))$, we need to work with a more general function, $h(u, b_H(u))$. More precisely, we would have to extend Theorem 3.7 to the function $h(u, x) = \sigma |c| \exp(\mu u + \sigma x)$. In return for which Theorem 3.10 follows. Another way consists in applying the Girsanov theorem first to the functional $\frac{1}{\sqrt{\varepsilon}} [\sqrt{\frac{\pi}{2}} \frac{\varepsilon^{(2-H)}}{\sigma_{2H}} \int_0^1 |\ddot{T}_\varepsilon(u)| \, \mathrm{d}u - \sigma \int_0^1 |X(u)| \, \mathrm{d}u] + d \int_0^1 |X(u)| \, \mathrm{d}u$, where $T_\varepsilon(\cdot) = \varphi_\varepsilon * X(\cdot)$, that is asymptotically equivalent to $F_\varepsilon$. Then, to show that this functional stably converges in case where $X(u) = c \exp(\sigma b_H(u))$ ($\mu(\cdot) \equiv 0$), we will need to extend Theorem 3.7 to the function $h(x) = \sigma |c| \exp(\sigma x)$. Convergence in Theorem 3.10 will then only be in distribution). □

**Remark.** *In Remark 3 following Theorem 3.9, using the Girsanov theorem we saw that for the third model we only obtain weak convergence for the case where $h(\cdot) \equiv 1$ and $k \geq 1$ and, a fortiori, for $k = 1$. If we apply Theorem 3.10 to this model under the true hypothesis $H_0$ and if we use Remark 2 following this theorem we note that this last convergence will take place stably. This is due to the fact that the computations are explicitly made in this part.*

### 3.3.2. About a variant on the last three models

Our techniques allow us also to consider the three following stochastic differential equations, for $t \geq 0$,

$$\mathrm{d}X_\varepsilon(t) = \sigma_\varepsilon(X_\varepsilon(t)) \, \mathrm{d}b_H(t) + \mu(X_\varepsilon(t)) \, \mathrm{d}t, \quad \text{with } X_\varepsilon(0) = c, \tag{36}$$

$X_\varepsilon(t) = c$, for $t < 0$ and we consider testing the hypothesis

$$H_0: \quad \sigma_\varepsilon(x) = \sigma \quad (\text{resp. } \sigma_\varepsilon(x) = \sigma x),$$

against the sequence of alternatives

$$H_\varepsilon: \quad \sigma_\varepsilon(x) = \sigma + \sqrt{\varepsilon}(d + F(\sqrt{\varepsilon}))x \quad (\text{resp. } \sigma_\varepsilon(x) = \sigma x + \sqrt{\varepsilon}(d + F(\sqrt{\varepsilon}))),$$

where $\sigma$, $d$ and $F(\cdot)$ are as in the previous section and $\mu(x) = \mu$ or $\mu(x) = \mu x$ for the first two models (resp. $\mu(x) = \mu x$ for the third one).

We will use the following result. Let the stochastic differential equation, for $t \geq 0$,

$$\mathrm{d}X(t) = (aX(t) + b) \, \mathrm{d}b_H(t) + \mu \, \mathrm{d}t$$

$$(\text{resp. } \mathrm{d}X(t) = (aX(t) + b) \, \mathrm{d}b_H(t) + \mu X(t) \, \mathrm{d}t),$$

with $X(0) = c$, and $a$, $b$ real constants such that $a \neq 0$. The solution is given by

$$X(t) = \frac{b}{a}\{\exp(ab_H(t)) - 1\} + \exp(ab_H(t))\left[\mu \int_0^t \exp(-ab_H(s))\,\mathrm{d}s + c\right]$$

$$\left(\text{resp. } X(t) = \frac{b}{a}\{\exp(ab_H(t)) - 1\} + \exp(\mu t + ab_H(t))\left[\frac{b\mu}{a}\int_0^t \exp(-\mu s)\{1 - \exp(-ab_H(s))\}\,\mathrm{d}s + c\right]\right).$$

By taking $a = \sqrt{\varepsilon}(d + F(\sqrt{\varepsilon}))$, $b = \sigma$ for the first two models (resp. $a = \sigma$ and $b = \sqrt{\varepsilon}(d + F(\sqrt{\varepsilon}))$ for the third one), for each model, there exists an unique process solution to the stochastic equation (36), let us say $X_\varepsilon(\cdot)$. Let us define as before the observed process $Y_\varepsilon(\cdot)$. We observe the following functionals

$$G_\varepsilon := \frac{1}{\sqrt{\varepsilon}}\left[\frac{\varepsilon^{2(2-H)}}{\sigma_{2H}^2}\int_0^1 \ddot{Y}_\varepsilon^2(u)\,\mathrm{d}u - \sigma^2\right]$$

$$\left(\text{resp. } G_\varepsilon := \frac{1}{\sqrt{\varepsilon}}\left[\frac{\varepsilon^{2(2-H)}}{\sigma_{2H}^2}\int_0^1 \ddot{Y}_\varepsilon^2(u)\,\mathrm{d}u - \sigma^2\int_0^1 Y_\varepsilon^2(u)\,\mathrm{d}u\right]\right).$$

Using a generalization of Theorems 3.6 and 3.7 we can prove the following theorem.

**Theorem 3.11.** *Let us suppose that $H$ is known with $1/2 < H < 1$, then $G_\varepsilon$ converges stably towards*

$$\sigma_{g_2}\sigma^2 N + 2d\sigma \int_0^1 X(u)\,\mathrm{d}u$$

$$\left(\text{resp. } G_\varepsilon \text{ converges stably towards } \sigma_{g_2}\sigma^2 \int_0^1 X^2(u)\,\mathrm{d}\widehat{W}(u) + 2d\sigma \int_0^1 X(u)\,\mathrm{d}u\right),$$

*where $X(\cdot)$ is the solution to (4) if $\mu(x) = \mu$ and to (5) if $\mu(x) = \mu x$ (resp. where $X(\cdot)$ is the solution to (6) and $\widehat{W}(\cdot)$ is a standard Brownian motion independent of $X(\cdot)$) and $g_2(\cdot)$ is defined by (10).*

**Remark.** *The two terms in the first sum are independent.*

**Proof of Theorem 3.11.** We just give an outline of the proof because this one follows the same lines as that one of Section 3.3.1. Indeed, we need to prove a lemma similar to Lemma 3.3 with attention paid to the case where $a = \sqrt{\varepsilon}(d + F(\sqrt{\varepsilon}))$ which tends to zero with $\varepsilon$. □

## 4. Final remark

The fourth model could be treated in the same way as the third one in Theorem 3.10. To obtain a similar result, it would be required to improve Theorem 3.7 for most general functions, say, $h(u, b_H(s), s \leq u) = \exp(\sigma b_H(u))|c + \mu \int_0^u \exp(-\sigma b_H(s))\,\mathrm{d}s|$.

The authors do not know at this moment if this kind of generalization can be done. Note that if this conjecture is true then the conjecture by Remark 4 after Theorem 3.9 would follow.

## Appendix

**Proof of Lemma 3.1.** We shall do the proof of this lemma for the third model. The other models could be treated in a similar way.

Using that $\int_{-\infty}^{\infty} \ddot{\varphi}(x)\,dx = 0$, one gets

$$\ddot{X}_\varepsilon(t) = \frac{c}{\varepsilon^2}\int_{-\infty}^{t/\varepsilon}\ddot{\varphi}(x)\exp[\sigma b_H(t-\varepsilon x)+\mu(t-\varepsilon x)]\,dx + \frac{c}{\varepsilon^2}\int_{t/\varepsilon}^{\infty}\ddot{\varphi}(x)\,dx$$

$$= \frac{c}{\varepsilon^2}\int_{-\infty}^{\infty}\ddot{\varphi}(x)\{\exp[\sigma b_H(t-\varepsilon x)+\mu(t-\varepsilon x)]-\exp[\sigma b_H(t)+\mu t]\}\,dx$$

$$+ \frac{c}{\varepsilon^2}\int_{t/\varepsilon}^{\infty}\ddot{\varphi}(x)\{1-\exp[\sigma b_H(t-\varepsilon x)+\mu(t-\varepsilon x)]\}\,dx.$$

Using once that $\int_{-\infty}^{\infty}\ddot{\varphi}(x)\,dx = \int_{-\infty}^{\infty}x\ddot{\varphi}(x)\,dx = 0$, and making a Taylor expansion of the exponential function one gets

$$\ddot{X}_\varepsilon(t) = \sigma X_\varepsilon(t)\ddot{b}_H^\varepsilon(t) + c_\varepsilon(t),$$

where

$$c_\varepsilon(t) = \sigma(X(t)-X_\varepsilon(t))\ddot{b}_H^\varepsilon(t)$$

$$+ \frac{c}{\varepsilon^2}\left[\int_{t/\varepsilon}^{\infty}\ddot{\varphi}(x)\{1-\exp[\sigma b_H(t-\varepsilon x)+\mu(t-\varepsilon x)]\}\,dx\right]\mathbf{1}_{\{0\le t\le\varepsilon\}}$$

$$+ \frac{c}{2\varepsilon^2}\int_{-\infty}^{\infty}\ddot{\varphi}(x)\{\sigma[b_H(t-\varepsilon x)-b_H(t)]-\mu\varepsilon x\}^2$$

$$\times \exp\{\sigma b_H(t)+\mu t+\theta[-\mu\varepsilon x+\sigma(b_H(t-\varepsilon x)-b_H(t))]\}\,dx$$

$$= (1)+(2)+(3),$$

with $0\le\theta<1$.

We are going to bound (1).

We know by [11] that $X(\cdot)$ is a $(H-\delta)$-Hölder continuous function on all compact included in $\mathbb{R}^+$ (see the remark following hypotheses (H1) and (H2) in Section 3.2), thus

$$|X_\varepsilon(t)-X(t)|\le \mathbf{C}(\omega)\varepsilon^{(H-\delta)}. \tag{A.1}$$

Furthermore, using the inequality (20), one gets

$$|(1)|\le \mathbf{C}(\omega)\varepsilon^{(2H-2-\delta)}.$$

Now to bound (2) and (3), we use the inequality $|1-\exp(x)|\le 2|x|$, for $|x|$ small enough and the modulus of continuity of $b_H(\cdot)$ (see (19)). We obtain that $|(2)|\le \mathbf{C}(\omega)\varepsilon^{(H-2-\delta)}$ and $|(3)|\le \mathbf{C}(\omega)\varepsilon^{(2H-2-\delta)}$.

Now to conclude the proof of this lemma, we have just to prove that $X_\varepsilon(\cdot)$ is uniformly bounded below on $[0,1]$.

By (A.1), for all $0\le t\le 1$ and any $\delta>0$,

$$||X_\varepsilon(t)|-|X(t)||\le |X_\varepsilon(t)-X(t)|\le \mathbf{C}(\omega)\varepsilon^{(H-\delta)},$$

and then,

$$|X_\varepsilon(t)|\ge |X(t)|-\mathbf{C}(\omega)\varepsilon^{(H-\delta)}.$$

Furthermore for all $t\ge 0$, $X(t) = c\exp(\sigma b_H(t)+\mu t)$ and then for all $t\in[0,1]$, $|X(t)|\ge |c|\exp(-a(\omega))>0$, where $a(\omega)=\sigma\sup_{t\in[0,1]}|b_H(t)(\omega)|+|\mu|$. Thus, we prove that $|X_\varepsilon(t)|\ge \frac{|c|}{2}\exp(-a(\omega))$ if $\varepsilon\le\varepsilon(\omega)$.

**Remark 3.** *For the fourth model, we also need to find a lower bound for $X_\varepsilon(\cdot)$ and it is the reason why we asked for $\mu$ and $c$ to have the same sign.* □

**Proof of Lemma 3.2.** Using that $X(t) = K(b_H(t))$ for $t \geq 0$ and that $X(t) = K(b_H(0)) = c$, for $t < 0$, the fact that $\int_{-\infty}^{\infty} \ddot{\varphi}(x)\,dx = 0$ and that $K(\cdot)$ is the solution to the ODE (31), one has

$$\varepsilon^2(\ddot{X}_\varepsilon(t) - \sigma(K(b_H^\varepsilon(t)))\ddot{b}_H^\varepsilon(t))$$
$$= \varepsilon^2(\ddot{X}_\varepsilon(t) - \dot{K}(b_H^\varepsilon(t))\ddot{b}_H^\varepsilon(t))$$
$$= \int_{-\infty}^{\infty} \ddot{\varphi}(v)[K(b_H(t - \varepsilon v)) - K(b_H^\varepsilon(t)) - \dot{K}(b_H^\varepsilon(t))\{b_H(t - \varepsilon v) - b_H^\varepsilon(t)\}]\,dv$$
$$+ \left(\int_{t/\varepsilon}^{\infty} \ddot{\varphi}(v)[K(b_H(0)) - K(b_H(t - \varepsilon v))]\,dv\right)\mathbf{1}_{\{0 \leq t \leq \varepsilon\}}.$$

Now, making a Taylor expansion for the function $K(\cdot)$, one gets

$$\varepsilon^2(\ddot{X}_\varepsilon(t) - \sigma(K(b_H^\varepsilon(t)))\ddot{b}_H^\varepsilon(t))$$
$$= \frac{1}{2}\int_{-\infty}^{\infty} \ddot{\varphi}(v)\ddot{K}(\theta_1)(b_H(t - \varepsilon v) - b_H^\varepsilon(t))^2\,dv$$
$$+ \left(\int_{t/\varepsilon}^{\infty} \ddot{\varphi}(v)\dot{K}(\theta_2)[b_H(0) - b_H(t - \varepsilon v)]\,dv\right)\mathbf{1}_{\{0 \leq t \leq \varepsilon\}},$$

where $\theta_1$ (resp. $\theta_2$) is a point between $b_H(t - \varepsilon v)$ and $b_H^\varepsilon(t)$ (resp. between $b_H(0)$ and $b_H(t - \varepsilon v)$). The modulus of continuity of $b_H(\cdot)$ (see (19)) yields the required inequality.

Similar computations could be done for $(X_\varepsilon(t) - K(b_H^\varepsilon(t)))$ and the lemma follows. □

**Proof of Lemma 3.3.** The proof concerning $\ddot{Y}_\varepsilon(t)$ is based on the proof of Lemma 3.1. It consists in bounding expressions (2) and (3) that appear in this lemma with $\sigma_\varepsilon$ taking the role of $\sigma$ using the fact that $\sigma_\varepsilon$ is bounded. Concerning the writing of $Y_\varepsilon(t) = \varphi_\varepsilon * X_\varepsilon(t)$, we use the expression of $X_\varepsilon(t)$ that is, $X_\varepsilon(t) = c\exp(\mu t + \sigma_\varepsilon b_H(t))$, the modulus of continuity of $b_H(\cdot)$ (see (19)) and the fact that $\sigma_\varepsilon$ is bounded. □